\documentclass[portuguese,a4paper,twoside,10pt]{article}
\usepackage[english]{babel}
\usepackage[T1]{fontenc}
\usepackage[latin1]{inputenc}
\usepackage{amssymb,amsmath,epsfig,amsfonts,euscript}
\usepackage{dsfont}
\usepackage[usenames,dvipsnames]{color}

\newtheorem{teo}{Theorem}[section]
\newtheorem{pro}[teo]{Proposition}

\newtheorem{lem}[teo]{Lema}

\newtheorem{ex}{Example}

\newcounter{note}[section]
\newcounter{example}[section]

\newcommand{\en}{\mathbb{N}}

\newcommand{\er}{\mathbb{R}}
\newcommand{\nd}{non degenerate}

\newcommand{\esp}{ }
\newcommand{\hs}{\hspace{3pt}}
\newcommand{\dst}{\displaystyle}
\newcommand{\dem}{{\bf Dem. }}
\newcommand{\fdem}{$\square$}

\newcommand{\nn}{\nonumber}
\newcommand{\mb}{\mathbf}

\newcommand{\titulo}[1]{\begin{center}\mbox{} \\ \noindent \textit{\textbf{\Large #1}}\\\vspace{0.5cm}\end{center}}

\renewcommand{\abstract}[1]{{\small \noindent \textbf{Abstract:} #1\\}}

\newcommand{\keywords}[1]{{\small \noindent \textbf{Keywords:} #1\\}}

\begin{document}

\begin{center}
\titulo{Extremal behavior of pMAX processes}
\end{center}

\vspace{0.5cm}

\textbf{Helena Ferreira} Department of Mathematics, University of
Beira
Interior, Covilhã, Portugal \\(helena.ferreira@ubi.pt)\\

\textbf{Marta Ferreira} Department of Mathematics, University of
Minho, Braga, Portugal \\(msferreira@math.uminho.pt)\\

\abstract{The well-known M4 processes of Smith and Weissman are very flexible models for asymptotically dependent multivariate data. Extended M4 of Heffernan \emph{et al.} allows to also account for asymptotic independence. In this paper we introduce a more general multivariate model comprising asymptotic dependence and independence, which has the extended M4 class as a particular case. We study properties of the proposed model. In particular, we compute
the multivariate extremal index, tail
dependence and extremal coefficients.}

\\


\section{Introduction}

Smith and Weissman (\cite{smith+weiss}, 1996) presented the so called
multivariate maxima of moving maxima (henceforth M4) process for the modeling of cross-sectional and serial extreme dependencies.
The class of M4 processes is very flexible for data which
exhibit asymptotic dependence. However, there are also many examples of data which are asymptotically independent (see Ledford and Tawn \cite{led+tawn1,led+tawn2} 1996/1997, among others), and hence cannot suitably be modeled by M4. Extended M4 processes (EM4) were latter considered in Heffernan \emph{et al.} (\cite{zangh+}, 2007) in order to also account for asymptotic independence.

The concept of tail independence between two random variables (r.v.'s), $X_1$ and $X_2$, with identical marginal distribution function (d.f.) $F$, was introduced in Sibuya (\cite{sib}, 1960). More precisely, they are said to be asymptotically independent if
 \begin{eqnarray}\label{lambda}
\lambda=\dst\lim_{x\to x_F}P(X_2 > x| X_1 > x)
\end{eqnarray}
is null and asymptotically dependent whenever $\lambda>0$, with $x_F=\sup\{x:F(x)<1\} $ being the upper end-point. Coefficient $\lambda$ quantifies the amount of dependence of the bivariate upper tails and is usually denoted \emph{tail dependence coefficient} (TDC). For instance, gaussian random pairs have $\lambda=0$, i.e., are asymptotically independent whereas $t$-distributed ones have $\lambda>0$ and thus asymptotically dependent.

In the asymptotic independent case, Ledford and Tawn (\cite{led+tawn1,led+tawn2} 1996/1997) proposed to model the null limit in (\ref{lambda}) by introducing a new coefficient ($\eta$) to rule the decay rate of the joint bivariate survival function
evaluated at the same large $x$. More precisely,
 \begin{eqnarray}\label{eta}
P(X_1 > x, X_2 > x)\sim L\Big(\frac{ 1}{P(X_1>x)}\Big)P(X_1 > x)^{1/\eta}, \textrm{ as $x \to x_F $,}
\end{eqnarray}
where $L$ is a slowly varying function, i.e. $L(tx)/L(x)\to 1$ as $x\to\infty$ for any
fixed $t > 0$, and $\eta\in (0,1]$ is a constant. The r.v.'s $X_1$ and $X_2$ are called positively associated when
$1/2 < \eta < 1$, nearly independent when $\eta= 1/2$ and negatively associated when $0 < \eta < 1/2$. Observe that they are asymptotically dependent if $\eta=1$ and $L(x)\not \to 0$, as $x\to\infty$, and asymptotically independent otherwise. Explicit formulas for $\eta$ can be seen in Heffernan (\cite{hef}, 2000)  for several known joint distributions.

The lag-$r$ ($r\in\en_0$) tail dependence coefficient for stationary $d$-dimensional sequences, $\{\mb{X}_n=(X_{n,1},...,X_{n,d})\}_{n\geq 1}$,
 with identical marginal distribution $F$, is naturally stated as (see Heffernan \emph{et al.} \cite{zangh+}, 2007 and references therein)
 \begin{eqnarray}\label{lambdar}
\lambda_{jj'}^{(r)}(\mb{X})=\lim_{x\to x_F}P(X_{1+r,j'} > x| X_{1,j} > x).
\end{eqnarray}
Analogously, the lag-$r$ ($r\in\en_0$) Ledford and Tawn  coefficient $\eta_{jj'}^{(r)}(\mb{X})$ is defined by
 \begin{eqnarray}\label{etar}
P(X_{1,j} > x, X_{1+r,j'} > x)\sim L\Big(\frac{ 1}{P(X_{1,j}>x)}\Big)P(X_{1,j} > x)^{1/\eta_{jj'}^{(r)}(\mb{X})}, \textrm{ as $x \to x_F $}.
\end{eqnarray}
We have that $\lambda_{jj'}(\mb{X})\equiv \lambda_{jj'}^{(0)}(\mb{X})$ is the TDC between the $j$th and the $j'$th components, $\lambda_{j}^{(r)}(\mb{X})\equiv \lambda_{jj}^{(r)}(\mb{X})$ is the lag-$r$ TDC within the $j$th sequence and $\lambda_{jj'}^{(r)}(\mb{X})$
is the lag-$r$ cross-sectional TDC between the $j$th and the $j'$th sequences. A similar deduction concerns the Ledford and Tawn coefficients, respectively, $\eta_{jj'}(\mb{X})$, $\eta_{j}^{(r)}(\mb{X})$ and $\eta_{jj'}^{(r)}(\mb{X})$.\\

A phenomenon also noticed in real data is that extreme events often tend to occur in clusters. The measure able to capture the clustered extremal dependence is the so called extremal index (Leadbetter \emph{et al.} \cite{lead+} 1983). We shall define the multivariate extremal index (Nandagopalan \cite{nand}, 1990) from which we can derive the univariate one of each marginal component. If
$\{\mathbf{\widehat{X}}_n=(\widehat{X}_{n,1},...,\widehat{X}_{n,d})\}_{n\geq
1}$ is an i.i.d.\hs sequence such that, for some sequences of
constants, $\{\mathbf{a}_n=(a_{n1}>0,...,a_{nd}>0)\}_{n\geq 1}$ and
$\{\mathbf{b}_n=(b_{n1},...,b_{nd})\}_{n\geq 1}$, the vector of
componentwise maxima
$\mathbf{\widehat{M}}_n=(\widehat{M}_{n\,1},...,\widehat{M}_{n\,d})$
satisfies
\begin{eqnarray}\label{mevlimit}
F_{\mb{\widehat{X}}}^n(\mathbf{a}_n\mathbf{x}+\mathbf{b}_n)\equiv P(\mathbf{\widehat{M}}_n\leq
\mathbf{a}_n\mathbf{x}+\mathbf{b}_n)
\dst\mathop{\longrightarrow}_{n\to\infty} H(x_1,...,x_d),
\end{eqnarray}
with $H$ a \nd\esp d.f.,  then $H$ is a multivariate extreme value d.f. (MEV) and we say that $F_{\mb{\widehat{X}}}$  belongs to the domain of attraction of $H$, in short  $F_{\mb{\widehat{X}}}\in\mathcal{D}(H)$. We recall that the convergence in (\ref{mevlimit})  can also be stated using the copula concept, i.e., $
C_{F_{\mb{\widehat{X}}}}(F_{\widehat{X}_1}(X_1),...,F_{\widehat{X}_d}(X_d))=F_{\mb{\widehat{X}}}(x_1,...,x_d)$, $(x_1,...,x_d)\in\er^d$. It holds $F_{\mb{\widehat{X}}}\in\mathcal{D}(H)$ if and only if $F_{\widehat{X}_j}\in\mathcal{D}(H_j)$, $j=1,...,d$, and
\begin{eqnarray}\nn
\begin{array}{lc}
C_{F_{\mb{\widehat{X}}}}^n(u_1^{1/n},...,u_d^{1/n})\dst\mathop{\longrightarrow}_{n\to\infty} C_H(u_1,...,u_d).
\end{array}
\end{eqnarray}
A stationary sequence $\{\mathbf{X}_n\}_{n\geq 1}$, having
common distribution $F_{\mathbf{X}}=F_{\mathbf{\widehat{X}}}$, has extremal index $\theta(\boldsymbol{\tau})\equiv
\theta(\tau_1,...,\tau_d)\in[0,1]$ when, for each
$\boldsymbol{\tau}=(\tau_1,...,\tau_d)\in\mathbb{R}_+^d$, there
exists
$\{\mathbf{u}_n^{(\boldsymbol{\tau})}={(u_{n,j}^{(\tau_j)},1\leq
j\leq d)}\}_{n\geq 1}$, satisfying 
\begin{eqnarray}\label{normlev}
n(1-F_{X_{n,j}}({u_{n,j}^{(\tau_j)}}))
\dst\mathop{\longrightarrow}_{n\to\infty}\tau_j,\,j=1,...,d,
\end{eqnarray}

\begin{eqnarray}\nn
P(\mathbf{\widehat{M}}_n\leq \mathbf{u}_n^{(\boldsymbol{\tau})})
\dst\mathop{\longrightarrow}_{n\to\infty}\gamma(\boldsymbol{\tau})
\textrm{ and }
P(\mathbf{M}_n\leq \mathbf{u}_n^{(\boldsymbol{\tau})})
\dst\mathop{\longrightarrow}_{n\to\infty}
\gamma(\boldsymbol{\tau})^{\theta(\boldsymbol{\tau})}.
\end{eqnarray}
Vectors $\mathbf{u}_n^{(\boldsymbol{\tau})}$ satisfying (\ref{normlev}) are usually denoted normalized levels.
Just as in one dimension, the extremal index is a key parameter
relating the extreme value properties of a stationary sequence
$\{\mathbf{X}_n\}_{n\geq 1}$ to those of the i.i.d.\hs associated
sequence  $\{\mathbf{\widehat{X}}_n\}_{n\geq 1}$.
If (\ref{mevlimit}) holds and $\{\mathbf{X}_n\}_{n\geq 1}$ has
multivariate extremal index $\theta(\boldsymbol{\tau})$, then
\begin{eqnarray}\nn
P(\mathbf{M}_n\leq \mathbf{a}_n\mathbf{x}+\mathbf{b}_n)
\dst\mathop{\longrightarrow}_{n\to\infty}
G(x_1,...,x_d),
\end{eqnarray}
and the MEV d.f.\hs in the limit satisfies:
\begin{eqnarray}\label{mevdepindep}
\begin{array}{c}
G(x_1,...,x_d)
=H(x_1,...,x_d)^{\theta(\tau_1(x_1),...,\tau_d(x_d))}\\
\end{array}
\end{eqnarray}
and
\begin{eqnarray}\nn
\begin{array}{c}
G_{j}(x_j)=H_j^{\theta_j}(x_j),
\end{array}
\end{eqnarray}
with
\begin{eqnarray}\nn
\begin{array}{ccc}
\tau_j(x_j)=-\log H_j(x_j),\,j=1,...,d,& \textrm{ and }& \theta_j=\dst\lim_{\stackrel{\tau_i\to
0}{i\not=j}}\theta(\tau_1,...,\tau_d).
\end{array}
\end{eqnarray}
The multivariate extremal index, although dependent of $\tau$,
satisfies the following property:
\begin{eqnarray}\label{hom}
\theta(c\tau_1,...,c\tau_d)=\theta(\tau_1,...,\tau_d),\,\forall c>0.
\end{eqnarray}
In a univariate context, $\theta_j$ is called the extremal index of the sequence $\{X_{n,j}\}_{n\geq 1}$.
The value of $1/\theta_j$ can be interpreted as the limiting mean number of exceedances of a
threshold per independent cluster as the threshold increases. A unit $ \theta_j$ means no serial clustering and is a form of asymptotic independence of extremes.

Here we propose the pMAX model constructed from a $d$-dimensional stationary sequence $\{\mb{X}_n\}_{n\geq 1}=\mb{X}$ with multivariate extremal index and a sequence of independent and identically distributed (i.i.d) random vectors $\{\mb{Z}_n\}_{n\geq 1}=\mb{Z}$, whose marginals are transformed through a positive exponent. The variation of the exponent values allows to obtain models ranging from serial asymptotic independence to asymptotic dependence. The pMAX sequence may also exhibit clustering of large values. In computing
the multivariate extremal index, tail
dependence coefficients and the extremal coefficient we shall encounter a very rich and wide class that includes the EM4 processes {(Section \ref{spmax}). We end with some brief notes on the estimation within the new model (Section \ref{sestim})}.

\section{The pMAX model}\label{spmax}

Consider $\{\mb{X}_n=(X_{n,1},...,X_{n,d})\}_{n\geq 1}$ a stationary sequence with multivariate extremal index $\theta^{\mb{X}}(\tau_1,...,\tau_d)$, $(\tau_1,...,\tau_d)\in\er_0^+$, and common distribution function $F_{\mb{{X}}}$ with unit Fréchet marginals  $F_{X_j}(x)=\exp(-1/x)$, $x > 0$, and in the domain of attraction of $G$, i.e., $C_{F_{\mb{{X}}}}^n(u_1^{1/n},...,u_d^{1/n})\dst\mathop{\longrightarrow}_{n\to\infty} C_G(u_1,...,u_d)$.  Consider an i.i.d. sequence of random vectors $\{\mb{Z}_n=(Z_{n,1},...,Z_{n,d})\}_{n\geq 1}$, also with unit Fréchet marginals $F_{Z_j}(x)=\exp(-1/x)$, $x > 0$,
and satisfying $F_{\mb{{Z}}}\in\mathcal{D}(H)$. Let $\boldsymbol{\alpha}=(\alpha_1,...,\alpha_d)$ be a vector of real positive constants.

The pMAX model is defined by
\begin{eqnarray}\label{pmax}
\begin{array}{lc}
\{\mb{Y}_n=(Y_{n,1},...,Y_{n,d})\}_{n\geq 1}=\{(X_{n,1}\vee Z_{n,1}^{1/\alpha_1},...,X_{n,d}\vee Z_{n,d}^{1/\alpha_d})\}_{n\geq 1},
\end{array}
\end{eqnarray}
with notation $a\vee b=\max(a,b)$. {The letter ``p" stands for the power transformation concerning $\{\mb{Z}_n\}_{n\geq 1}$}. The EM4 model (Heffernan \emph{et al.}, \cite{zangh+} 2007) is obtained by considering $\{\mb{X}_n\}_{n\geq 1}$ an M4 process, $\alpha_j=\alpha$, $j=1,...,d$ and $Z_{n,1},...,Z_{n,d}$ independent.

At this point, we state some notation used along the paper. We will denote $U$  the set of indexes $j$ in $D=\{1,...,d\}$ for which $\alpha_j\geq 1$. For any $A\subset D$, the vector $(x_1,...,x_d)_A$ denotes the sub-vector of $(x_1,...,x_d)$ with indexes in $A$ and, for any d.f. $F$, $F_A$ denotes the marginal d.f. of $F$ for sub-vectors with indexes in $A$. Also $\mathbf{{M}}_n$ and $\mathbf{\widehat{M}}_n$ will denote the componentwise maxima of, respectively, the pMAX sequence $\{\mb{Y}_n\}_{n\geq 1}$ and  the corresponding i.i.d.\hs
sequence $\{\mb{\widehat{Y}}_n\}_{n\geq 1}$. {We will also use notation $a\wedge b=\min(a,b)$}.

\subsection{The extremal index}

In this section we compute the marginal and the multivariate extremal indexes of a pMAX process. We will find the interesting and novel feature (when comparing with the EM4) that pMAX may have multivariate extremal index not constant equal to one while some marginal extremal indexes may be unit.

First we need to compute the normalized levels form of pMAX.
\begin{pro}\label{pnormlev}
A sequence of normalized levels $\{\mathbf{u}_n^{(\boldsymbol{\tau})}={(u_{n,j}^{(\tau_j)},1\leq
j\leq d)}\}_{n\geq 1}$ for a pMAX process is such that, for each $j\in D$,
\begin{eqnarray}\nn
u_{n,j}^{(\tau_j)}=\left\{
\begin{array}{ll}
\big(\frac{n}{\tau_j}\big)^{1/\alpha_j} &, \alpha_j<1\vspace{0.25cm}\\
\frac{n}{\tau_j} &, \alpha_j\geq 1\,.
\end{array}\right.
\end{eqnarray}
\end{pro}
\dem
First observe that
\begin{eqnarray}\nn
F_{Y_{n,j}}(x)=P(X_{n,j}\leq x)P(Z_{n,j}\leq x^{\alpha_j})=e^{-x^{-1}}e^{-x^{-\alpha_j}}.
\end{eqnarray}
If $\alpha_j<1$, then
$$
\lim_{n\to\infty}F_{Y_{n,j}}^n\Big(\Big(\frac{n}{\tau_j}\Big)^{1/\alpha_j}\Big)
=\lim_{n\to\infty}\exp\Big(-\frac{\tau_j^{1/\alpha_j}}{n^{1/\alpha_j-1}}-\tau_j\Big)
= \exp(-\tau_j),
$$
whereas for $\alpha_j\geq 1$, we have
$$
\lim_{n\to\infty}F_{Y_{n,j}}^n\Big(\frac{n}{\tau_j}\Big)
=\lim_{n\to\infty}\exp\Big(-\tau_j-\frac{\tau_j^{\alpha_j}}{n^{\alpha_j-1}}\Big)
= \exp(-\tau_j).\,\,\square
$$

Now we compute the marginal extremal indexes of a pMAX process.

\begin{pro}\label{ppmaxeiuniv}
For each $j\in D$, if $\alpha_j\geq 1$ then $\{Y_{n,j}\}_{n\geq 1}$ has extremal index coincident with the one of $\{X_{n,j}\}_{n\geq 1}$, otherwise it will be unit.
\end{pro}
\dem
For the case $\alpha_j<1$, we have normalized levels $u_{n,j}^{(\tau_j)}=(n/\tau_j)^{1/\alpha_j}$ (Proposition \ref{pnormlev}) and hence
$$
\begin{array}{c}
P(M_{n,j}\leq u_{n,j}^{(\tau_j)})=P(\bigvee_{i=1}^n X_{i,j}\leq (n/\tau_j)^{1/\alpha_j})P(\bigvee_{i=1}^n Z_{i,j}\leq n/\tau_j)
\end{array}
$$
Observe that
$$
\begin{array}{c}
\lim_{n\to\infty}F_{X_{i,j}}^n\Big(\Big(\frac{n}{\tau_j}\Big)^{1/\alpha_j}\Big)
=\lim_{n\to\infty}\exp\Big(-\frac{\tau_j^{1/\alpha_j}}{n^{1/\alpha_j-1}}\Big)
= 1
\end{array}
$$
and thus
$$
\begin{array}{c}
\lim_{n\to\infty}P(\bigvee_{i=1}^n X_{i,j}\leq (n/\tau_j)^{1/\alpha_j})=e^{-\theta_j\times 0}=1.
\end{array}
$$
Therefore,
\begin{eqnarray}\nn
\begin{array}{c}
\lim_{n\to\infty}P(M_{n,j}\leq u_{n,j}^{(\tau_j)})=\lim_{n\to\infty}P(\bigvee_{i=1}^n Z_{i,j}\leq n/\tau_j)=e^{-\tau_j}.
\end{array}
\end{eqnarray}

If $\alpha_j\geq 1$, the normalized levels $u_{n,j}^{(\tau_j)}=n/\tau_j$ (Proposition \ref{pnormlev}) and hence
\begin{eqnarray}\nn
\begin{array}{c}
\lim_{n\to\infty}P(M_{n,j}\leq u_{n,j}^{(\tau_j)})=\lim_{n\to\infty}P(\bigvee_{i=1}^n X_{i,j}\leq (n/\tau_j)^{1/\alpha_j})P(\bigvee_{i=1}^n Z_{i,j}\leq n/\tau_j)=e^{-\theta_j\tau_j}.\,\, \square\\\\
\end{array}
\end{eqnarray}

In what follows we adopt the convention that if $U=\emptyset$ or $D-U=\emptyset$ the quantities involving the respective sub-distributions are considered null.

\begin{pro}\label{ppmaxei}
The multivariate extremal index of $\{\mb{Y}_n\}_{n\geq 1}$ is given by
\begin{eqnarray}\label{pmaxei}
\begin{array}{c}
\dst\theta^{\mb{Y}}(\tau_1,...,\tau_d)
=\frac{\theta^{\mb{X}}_U(\tau_1,...,\tau_d)_U\log C_{G_U}(e^{-\tau_1},...,e^{-\tau_d})_U+\log C_{H_{D-U}}(e^{-\tau_1},...,e^{-\tau_d})_{D-U}}{\log C_{G_U}(e^{-\tau_1},...,e^{-\tau_d})_U+\log C_{H_{D-U}}(e^{-\tau_1},...,e^{-\tau_d})_{D-U}},
\end{array}
\end{eqnarray}
where
\begin{eqnarray}\nn
\begin{array}{c}
\theta^{\mb{X}}_U(\tau_1,...,\tau_d)_U=\dst\lim_{\tau_j\downarrow 0,j\in D-U}\theta^{\mb{X}}(\tau_1,...,\tau_d).
\end{array}
\end{eqnarray}
\end{pro}
\dem
We have, successively,
\begin{eqnarray}\nn
\begin{array}{rl}
&P(\mb{M}_n\leq \mb{u}_n^{(\boldsymbol{\tau})})
=P(M_{n,1}\leq {u}_{n,1}^{({\tau_1})},...,M_{n,d}\leq {u}_{n,d}^{({\tau_d})})\vspace{0.35cm}\\
=& P(\{\bigvee_{i=1}^nX_{i,j}\leq \frac{n}{\tau_j},j\in U\},
\{\bigvee_{i=1}^nX_{i,j}\leq \big(\frac{n}{\tau_j}\big)^{1/\alpha_j},j\in D-U\})\times\vspace{0.15cm}\\
&
P(\{\bigvee_{i=1}^nZ_{i,j}\leq \big(\frac{n}{\tau_j}\big)^{\alpha_j} ,j\in U\},
\{\bigvee_{i=1}^nZ_{i,j}\leq \frac{n}{\tau_j} ,j\in D-U\})\vspace{0.35cm}\\
=&\Big( P(\{\bigvee_{i=1}^nX_{i,j}\leq \frac{n}{\tau_j},j\in U\})-P(\{\bigvee_{i=1}^nX_{i,j}\leq \frac{n}{\tau_j},j\in U\},
\bigcup_{j\in D-U}\{\bigvee_{i=1}^nX_{i,j}> \big(\frac{n}{\tau_j}\big)^{1/\alpha_j}\})\Big)\times\vspace{0.15cm}\\
&
\Big(P(\{\bigvee_{i=1}^nZ_{i,j}\leq \frac{n}{\tau_j} ,j\in D-U\})-P(\{\bigvee_{i=1}^nZ_{i,j}\leq \frac{n}{\tau_j} ,j\in D-U\},\bigcup_{j\in U}\{\bigvee_{i=1}^nZ_{i,j}> \big(\frac{n}{\tau_j}\big)^{\alpha_j}\})\Big)
\end{array}
\end{eqnarray}
Since the second terms in each of the two factors converge to zero, we have
\begin{eqnarray}\nn
\begin{array}{rl}
{\dst \lim_{n\to\infty}}P(\mb{M}_n\leq \mb{u}_n^{(\boldsymbol{\tau})})
=&{\dst \lim_{n\to\infty}}P(\bigvee_{i=1}^nX_{i,j}\leq \frac{n}{\tau_j},j\in U)P(\bigvee_{i=1}^nZ_{i,j}\leq \frac{n}{\tau_j} ,j\in D-U)\vspace{0.35cm}\\
=&C_{G_U}^{\theta^{\mb{X}}_U(\tau_1,...,\tau_d)_U}(e^{-\tau_1},...,e^{-\tau_d})_U\,\,C_{H_{D-U}}(e^{-\tau_1},...,e^{-\tau_d})_{D-U}.
\end{array}
\end{eqnarray}
On the other hand,
\begin{eqnarray}\label{mevindep}
\begin{array}{rl}
{\dst \lim_{n\to\infty}}P(\mb{\widehat{M}}_n\leq \mb{u}_n^{(\boldsymbol{\tau})})
=&{\dst\lim_{n\to\infty}}P^n(Y_{n,j}\leq{u}_{n,j}^{({\tau_j})},j\in D)\vspace{0.35cm}\\
=&{\dst\lim_{n\to\infty}}P^n(X_{n,j}\leq \frac{n}{\tau_j},j\in U)P^n(Z_{n,j}\leq \frac{n}{\tau_j},j\in D-U)\vspace{0.35cm}\\
=&C_{G_U}(e^{-\tau_1},...,e^{-\tau_d})_U\,\,C_{H_{D-U}}(e^{-\tau_1},...,e^{-\tau_d})_{D-U}
\end{array}
\end{eqnarray}

Now the result follows from relation (\ref{mevdepindep}).\,\fdem\\\\

Observe that, from expression (\ref{pmaxei}), we obtain the result of Proposition \ref{ppmaxeiuniv} for the univariate $\theta_j^{{Y}}$, $j\in \{1,...,d\}$.

We can also conclude that, if $\alpha_j\geq 1$ for all $j\in D$ then the multivariate extremal index of $\{\mb{Y}_n\}_{n\geq 1}$  coincides with the multivariate extremal index of $\{\mb{X}_n\}_{n\geq 1}$. If $\alpha_j< 1$ for all $j\in D$ then the multivariate extremal index of $\{\mb{Y}_n\}_{n\geq 1}$  is unit. In the other cases, the multivariate extremal index of $\{\mb{Y}_n\}_{n\geq 1}$ depends on the multivariate extremal index of the sequence $\{(\mb{Y}_{U})_n\}_{n\geq 1}$ of the sub-vectors with components having indexes in $U$.

We shall see some particular cases:
\begin{itemize}
\item[(i)] If $C_G(u_1,...,u_d)=C_H(u_1,...,u_d)=\bigwedge_{j=1}^du_j$, then
 \begin{eqnarray}\nn
\begin{array}{c}
\dst\theta^{\mb{Y}}(\tau_1,...,\tau_d)
=\frac{\theta^{\mb{X}}_U(\tau_1,...,\tau_d)_U\bigvee_{j\in U}\tau_j+\bigvee_{j\in D-U}\tau_j}{\bigvee_{j\in U}\tau_j+\bigvee_{j\in D-U}\tau_j}.
\end{array}
\end{eqnarray}
\item[(ii)] If $C_G(u_1,...,u_d)=\bigwedge_{j=1}^du_j$ and $C_H(u_1,...,u_d)=\prod_{j=1}^du_j$, then
\begin{eqnarray}\nn
\begin{array}{c}
\dst\theta^{\mb{Y}}(\tau_1,...,\tau_d)
=\frac{\theta^{\mb{X}}_U(\tau_1,...,\tau_d)_U\bigvee_{j\in U}\tau_j+\sum_{j\in D-U}\tau_j}{\bigvee_{j\in U}\tau_j+\sum_{j\in D-U}\tau_j}.
\end{array}
\end{eqnarray}
\item[(iii)] If $C_G(u_1,...,u_d)=C_H(u_1,...,u_d)=\exp(-(\sum_{i=1}^d(-\log u_j)^{1/\beta})^{\beta})$, for some $0<\beta<1$, then
\begin{eqnarray}\nn
\begin{array}{c}
\dst\theta^{\mb{Y}}(\tau_1,...,\tau_d)
=\frac{\theta^{\mb{X}}_U(\tau_1,...,\tau_d)_U\big(\sum_{j\in U}\tau_j^{1/\beta}\big)^{\beta}+\big(\sum_{j\in D-U}\tau_j^{1/\beta}\big)^{\beta}}{\big(\sum_{j\in U}\tau_j^{1/\beta}\big)^{\beta}+\big(\sum_{j\in D-U}\tau_j^{1/\beta}\big)^{\beta}}.
\end{array}
\end{eqnarray}
\item[(iv)] If $\alpha_j=\alpha$, $j\in D$, then $\theta^{\mb{Y}}(\tau_1,...,\tau_d)=1$ whenever $\alpha<1$  and $\theta^{\mb{Y}}(\tau_1,...,\tau_d)=\theta^{\mb{X}}(\tau_1,...,\tau_d)$ if $\alpha\geq 1$. This corresponds to the case of EM4 model (see Heffernan \emph{et al.}, \cite{zangh+}, 2007).
\end{itemize}

\subsection{Extremal coefficients of the limiting MEV}

When computing the multivariate extremal index of the pMAX, we have already seen that if $F_{\mb{X}}\in\mathcal{D}(G)$ and $F_{\mb{Z}}\in\mathcal{D}(H)$ then
\begin{eqnarray}\nn
\begin{array}{c}
P(\mb{M}_n\leq \mb{u}_n^{(\boldsymbol{\tau})})
\dst\mathop{\longrightarrow}_{n\to\infty}
C_{G_U}^{\theta^{\mb{X}}_U(\tau_1,...,\tau_d)_U}(e^{-\tau_1},...,e^{-\tau_d})_U\,\,C_{H_{D-U}}(e^{-\tau_1},...,e^{-\tau_d})_{D-U}
\end{array}
\end{eqnarray}
Therefore, for $(x_1,...,x_d)\in\er_+^d$,
\begin{eqnarray}\nn
\begin{array}{rl}
&P(\{n^{-1}M_{n,j}\leq x_j, j\in U\},\{n^{-1/\alpha_j}M_{n,j}\leq x_j, j\in D-U\})\vspace{0.35cm}\\
\dst\mathop{\longrightarrow}_{n\to\infty}&
{G_U}^{\theta^{\mb{X}}_U(1/x_1,...,1/x_d)_U}(x_1,...,x_d)_U{H_{D-U}}(x_1^{\alpha_1},...,x_d^{\alpha_d})_{D-U}
\end{array}
\end{eqnarray}
We shall denote
\begin{eqnarray}\label{mevpmax}
\begin{array}{c}
V(x_1,...,x_d)={G_U}^{\theta^{\mb{X}}_U(1/x_1,...,1/x_d)_U}(x_1,...,x_d)_U{H_{D-U}}(x_1^{\alpha_1},...,x_d^{\alpha_d})_{D-U}.
\end{array}
\end{eqnarray}
We have also seen in (\ref{mevindep}) that $F_{\mb{Y}}\in\mathcal{D}(\widehat{V})$ where we have
\begin{eqnarray}\label{mevpmaxind}
\begin{array}{c}
\widehat{V}(x_1,...,x_d)=
{G_U}(x_1,...,x_d)_U{H_{D-U}}(x_1^{\alpha_1},...,x_d^{\alpha_d})_{D-U}.
\end{array}
\end{eqnarray}

The MEV distributions $\widehat{V}$ and $V$ do not have identically distributed marginals. Indeed, based on the proof of Proposition \ref{ppmaxeiuniv}, we may see that, for $\alpha_j<1$,
\begin{eqnarray}\nn
\begin{array}{c}
P(n^{-1/\alpha_j}M_{n,j}\leq x_j)\dst\mathop{\longrightarrow}_{n\to\infty}e^{-x_j^{-\alpha_j}}
\end{array}
\end{eqnarray}
and for $\alpha_j\geq 1$,
\begin{eqnarray}\nn
\begin{array}{c}
P(n^{-1}M_{n,j}\leq x_j)\dst\mathop{\longrightarrow}_{n\to\infty}e^{-\theta_jx_j^{-1}}.
\end{array}
\end{eqnarray}
Observe that we can now obtain these univariate extreme limiting distributions from the MEV ones, $V$ and $\widehat{V}$. \\

Since MEV distributions $G$ and $H$ have each identically distributed marginals, the definition of the extremal coefficient $\epsilon$ of Tiago de Oliveira (\cite{tiago}, 1962/63) and Smith (\cite{smith}, 1990) is pacific within this case. More precisely, for all $x>0$, we have
\begin{eqnarray}\nn
G(x,...,x)=(e^{-x^{-1}})^{\epsilon_G}\textrm{ and }
H(x,...,x)=(e^{-x^{-1}})^{\epsilon_H}
\end{eqnarray}
and thus
\begin{eqnarray}\nn
\epsilon_G= -\log G(1,...,1)\textrm{ and }
\epsilon_H= -\log H(1,...,1).
\end{eqnarray}
When the marginals of a MEV are not equally distributed, the previous definition needs some modification. A natural way is to consider an analogous definition based on MEV copulas. More precisely, the extremal coefficient of a MEV distribution $V$ is the constant $\epsilon_V$ such that
\begin{eqnarray}\label{eccop}
C_V(u,...,u)= u^{\epsilon_V},\, u\in[0,1].
\end{eqnarray}
If $V$ has identically distributed marginals, both definitions are equivalent.\\

We shall calculate the extremal coefficients of $\widehat{V}$ and $V$ according to definition in (\ref{eccop}).

\begin{pro}\label{ppmaxec}
The extremal coefficients of the MEV distributions, $V$ in (\ref{mevpmax}) and $\widehat{V}$ in (\ref{mevpmaxind}), are given by, respectively,
\begin{eqnarray}\nn
\epsilon_{\widehat{V}}=\epsilon_{G_U}+\epsilon_{H_{D-U}}.
\end{eqnarray}
and
\begin{eqnarray}\nn
\begin{array}{c}
\epsilon_{{V}}=\theta_U^{\mb{X}}(\frac{1}{\theta_1},...,\frac{1}{\theta_d})_U\Big(-\log C_{G_U}\big(e^{-1/\theta_1},...,e^{-1/\theta_d}\big)_U\Big)+\epsilon_{H_{D-U}}.
\end{array}
\end{eqnarray}
\end{pro}
\dem
Considering the variables change, $u_j=V(x_j)$ or $u_j=\widehat{V}_j(x_j)$, we obtain
\begin{eqnarray}\nn
C_{\widehat{V}}(u_1,...,u_d)= C_{G_U}(u_1,...,u_d)_U\,\,C_{H_{D-U}}(u_1,...,u_d)_{D-U},
\end{eqnarray}
which leads immediately to the first assertion, and
\begin{eqnarray}\nn
C_V(u_1,...,u_d)= C_{G_U}^{\theta^{\mb{X}}_U(-\log u_1^{1/\theta_1},...,-\log u_d^{1/\theta_d})_U}(u_1^{1/\theta_1},...,u_d^{1/\theta_d})_U\,\,C_{H_{D-U}}(u_1,...,u_d)_{D-U}
.
\end{eqnarray}
By the multivariate extremal index homogeneity property in (\ref{hom}), we have
\begin{eqnarray}\nn
C_{V}(u,...,u)= C_{G_U}^{\theta_U^{\mb{X}}(\frac{1}{\theta_1},...,\frac{1}{\theta_d})_U}(u^{1/\theta_1},...,u^{1/\theta_d})_U\,\,C_{H_{D-U}}(u,...,u)_{D-U},
\end{eqnarray}
If we use MEV's spectral representation of $G$ (Haan and Resnick, \cite{haan+res} 1977), then
\begin{eqnarray}\nn
\begin{array}{rl}
\dst C_{G_U}(u^{1/\theta_1},...,u^{1/\theta_d})_U
=& \dst\exp\Big(-\int_{S_d}\bigvee_{j\in U}\frac{w_j}{\theta_j}\,(-\log u)\,dW_G(w_1,...,w_2)\Big)\vspace{0.35cm}\\
=& u^{{\dst\int_{S_d}\bigvee_{j\in U}}\frac{w_j}{\theta_j}\,\,dW_G(w_1,...,w_2)},
\end{array}
\end{eqnarray}
where $W_G$ is a finite measure over the unit sphere $S_d$ of $\er_+^d$. Now just observe that
\begin{eqnarray}\nn
\begin{array}{c}
\epsilon_{{V}}=\theta_U^{\mb{X}}(\frac{1}{\theta_1},...,\frac{1}{\theta_d})_U\dst\int_{S_d}\bigvee_{j\in U}\frac{w_j}{\theta_j}\,\,dW_G(w_1,...,w_2)+\epsilon_{H_{D-U}}.\,\,\square\\\\
\end{array}
\end{eqnarray}

Let us apply the result in the case where $\{\mb{X}_n\}_{n\geq 1}$ is an M4 process. We have
$$
C_G(u_1,...,u_d)=\prod_{l=1}^\infty\prod_{k=-\infty}^\infty \bigwedge_{j=1}^du_j^{a_{lkj}},
$$
where $a_{lkj}$, $l\geq 1$, $-\infty<k<\infty$, $j\in D$, are the model constants and
$$
\theta^{\mb{X}}(\tau_1,...,\tau_d)=
\frac{\sum_{l=1}^\infty\bigvee_{j=1}^d\big(\tau_j
\bigvee_{k=-\infty}^\infty a_{lkj}\big)}
{\sum_{l=1}^\infty\sum_{k=-\infty}^\infty\bigvee_{j=1}^d(\tau_j\,
 a_{lkj})}.
$$
Therefore, any choice of $H$ will lead us to
$$
\begin{array}{ccc}
\dst\epsilon_{\widehat{V}}=\sum_{l=1}^\infty\sum_{k=-\infty}^\infty\bigvee_{j\in U} a_{lkj}+\epsilon_{H_{D-U}}
&\textrm{
and
}
&\dst\epsilon_{{V}}=\sum_{l=1}^\infty\bigvee_{k=-\infty}^\infty\bigvee_{j\in U} \frac{a_{lkj}}{\theta_j}+\epsilon_{H_{D-U}},
\end{array}
$$
where
$$
\theta_j=\sum_{l=1}^\infty\bigvee_{k=-\infty}^\infty a_{lkj},\,\, j=1,...,d.
$$
If in particular we choose $H$ as in Heffernan \emph{et al.} (\cite{zangh+}, 2007), we have $\epsilon_{\widehat{V}}=\epsilon_{{V}}=\epsilon_{H}=d$ if $\alpha<1$ and, for $\alpha\geq 1$,
$$
\begin{array}{ccc}
\dst\epsilon_{\widehat{V}}=\sum_{l=1}^\infty\sum_{k=-\infty}^\infty\bigvee_{j\in U} a_{lkj}
&\textrm{
and
}
&\dst\epsilon_{{V}}=\sum_{l=1}^\infty\bigvee_{k=-\infty}^\infty\bigvee_{j\in U} \frac{a_{lkj}}{\theta_j}.
\end{array}
$$

\subsection{Tail (in)dependence coefficients}

In this section we characterize the tail (in)dependence of a pMAX model $\{\mb{Y}_n\}_{n\geq 1}$ in (\ref{pmax}). We will compute the lag-$r$ TDC, $\lambda_{jj'}^{(r)}(\mb{Y})$ and the lag-$r$ Ledford and Tawn coefficient, $\eta_{jj'}^{(r)}(\mb{Y})$, by assuming that the respective exist for the underlying process $\{\mb{X}_n\}_{n\geq 1}$.\\

For calculations, we shall need the pMAX 
bivariate d.f.
\begin{eqnarray}\nn
\begin{array}{rl}
&\dst P(Y_{1,j}\leq x_j,Y_{1+r,j'}\leq x_{j'})=P(X_{1,j}\leq x_j,Z_{1,j}\leq x_j^{\alpha_j},X_{1+r,j'}\leq x_{j'},Z_{1+r,j'}\leq x_{j'}^{\alpha_{j'}})\vspace{0.35cm}\\
=& \dst P(X_{1,j}\leq x_{j},X_{1+r,j'}\leq x_{j'})P(Z_{1,j}\leq x_{j}^{\alpha_{j}})P(Z_{1+r,j'}\leq x_{j'}^{\alpha_{j'}})\vspace{0.35cm}\\
=& \dst P(X_{1,j}\leq x_{j},X_{1+r,j'}\leq x_{j'})e^{-x_{j}^{-\alpha_j}}e^{-x_{j'}^{-\alpha_{j'}}}.\\\\
\end{array}
\end{eqnarray}

\begin{pro}\label{ppmaxtdc}
If the pMAX process $\{\mb{Y}_n\}_{n\geq 1}$ in (\ref{pmax}) is such that $\lambda_{jj'}^{(r)}(\mb{X})$ exists for the underlying process $\{\mb{X}_n\}_{n\geq 1}$, then
\begin{eqnarray}\nn
\lambda_{jj'}^{(r)}(\mb{Y})=\left\{
\begin{array}{ll}
0 & ,\, \alpha_j<1\vspace{0.15cm}\\
\frac{1}{2}\lambda_{jj'}^{(r)}(\mb{X})& ,\, \alpha_j=1\vspace{0.15cm}\\
\lambda_{jj'}^{(r)}(\mb{X})& ,\, \alpha_j>1\vspace{0.15cm}.
\end{array}\right.
\end{eqnarray}
\end{pro}
\dem
We have
\begin{eqnarray}\nn
\begin{array}{rl}
&\dst P(Y_{1+r,j'}>x|Y_{1,j}>x)=\frac{1-P(Y_{1,j}\leq x)-P(Y_{1+r,j'}\leq x)+P(Y_{1,j}\leq x,Y_{1+r,j'}\leq x)}{1-P(Y_{1,j}\leq x)}\vspace{0.35cm}\\
=&\dst 1-\frac{P(X_{1+r,j'}\leq x)e^{-x^{-\alpha_{j'}}}+\big(P(X_{1+r,j'}\leq x)-P(X_{1,j}> x,X_{1+r,j'}\leq x)\big)e^{-x^{-\alpha_j}-x^{-\alpha_{j'}}}}{1-P(X_{1,j}\leq x)e^{-x^{-\alpha_j}}} \vspace{0.35cm}\\
=& \dst 1-e^{-x^{-1}-x^{-\alpha_{j'}}}\frac{1-e^{-x^{-\alpha_{j}}}}{1-e^{-x^{-1}-x^{-\alpha_j}}}+P(X_{1+r,j'}\leq x|X_{1,j}> x)e^{-x^{-\alpha_j}-x^{-\alpha_{j'}}}\frac{1-e^{-x^{-1}}}{1-e^{-x^{-1}-x^{-\alpha_j}}}.
\end{array}
\end{eqnarray}
Now just observe that
\begin{eqnarray}\nn
\lim_{x\to\infty}\frac{1-e^{-x^{-\alpha_{j}}}}{1-e^{-x^{-1}-x^{-\alpha_j}}}
=\left\{\begin{array}{ll}
1&,\,\alpha_j<1\\
1/2&,\,\alpha_j=1\\
0&,\,\alpha_j>1
\end{array}\right.
\end{eqnarray}
and
\begin{eqnarray}\nn
\lim_{x\to\infty}\frac{1-e^{-1}}{1-e^{-x^{-1}-x^{-\alpha_j}}}
=\left\{\begin{array}{ll}
0&,\,\alpha_j<1\\
1/2&,\,\alpha_j=1\\
1&,\,\alpha_j>1.\,\,\square
\end{array}\right.
\end{eqnarray}
\vspace{0.15cm}


In the particular case of the EM4 in Heffernan \emph{et al.} (\cite{zangh+}, 2007), we have $\alpha_j=\alpha$, $j=1,...,d$, leading to
\begin{eqnarray}\nn
\lambda_{jj'}^{(r)}(\mb{Y})=\left\{
\begin{array}{ll}
0 & ,\, \alpha<1\vspace{0.15cm}\\
\frac{1}{2}(2-\delta)& ,\, \alpha=1\vspace{0.15cm}\\
2-\delta& ,\, \alpha>1\vspace{0.15cm}.
\end{array}\right.
\end{eqnarray}
where $\delta=\sum_{l=1}^\infty\sum_{k=-\infty}^{\infty}a_{lkj}^{-1}\vee a_{lkj'}^{-1}$ (see Result 9 in Heffernan \emph{et al.}, \cite{zangh+} 2007; observe that a small correction is needed for the case $\alpha=1$).

\begin{pro}\label{ppmaxeta}
If the pMAX process $\{\mb{Y}_n\}_{n\geq 1}$ in (\ref{pmax}) is such that $\eta_{jj'}^{(r)}(\mb{X})$ exists for the underlying process $\{\mb{X}_n\}_{n\geq 1}$, then
\begin{eqnarray}\nn
\eta_{jj'}^{(r)}(\mb{Y})=\left\{
\begin{array}{ll}
\dst \max\Big(\frac{\alpha_j}{\alpha_j+\min(1,\alpha_{j'})},\alpha_j\eta_{jj'}^{(r)}(\mb{X})\Big)&,\,\alpha_j< 1\vspace{0.15cm}\\
\dst \max\Big(\frac{1}{1+\alpha_j},\frac{1}{1+\alpha_{j'}},\eta_{jj'}^{(r)}(\mb{X})\Big)&,\,\alpha_j\geq 1\,.
\end{array}\right.\vspace{0.35cm}
\end{eqnarray}
\end{pro}
\dem
By hypothesis,
\begin{eqnarray}\label{etax}
P(X_{1,j}> x,X_{1+r,j'}> x)\sim L(x)x^{-1/\eta_{jj'}^{(r)}(\mb{X})}
\end{eqnarray}
for some slowly varying function $L$. \\

We have
\begin{eqnarray}\nn
\begin{array}{c}
\dst P(Y_{1,j}>x,Y_{1+r,j'}>x)=
\dst 1-e^{-x^{-1}-x^{-\alpha_j}}-e^{-x^{-1}-x^{-\alpha_{j'}}}+P(X_{1,j}\leq x,X_{1+r,j'}\leq x)e^{-x^{-\alpha_j}-x^{-\alpha_{j'}}}
\end{array}
\end{eqnarray}
where
\begin{eqnarray}\nn
\begin{array}{rl}
&\dst P(X_{1,j}\leq x,X_{1+r,j'}\leq x)=2e^{-x^{-1}}-1+P(X_{1,j}> x,X_{1+r,j'}> x).
\end{array}
\end{eqnarray}

Observe now that, for any positive real constants $a,b$ and large enough $x$,
\begin{eqnarray}\nn
\begin{array}{rl}
\dst e^{-x^{-a}-x^{-b}}\sim 1-\frac{1}{x^a}-\frac{1}{x^b}+\frac{1}{2x^{2a}}+\frac{1}{2x^{2b}}+\frac{1}{2x^{a+b}}.
\end{array}
\end{eqnarray}
Therefore, after some calculations, we obtain
\begin{eqnarray}\nn
\begin{array}{c}
\dst \frac{P(Y_{1,j}>x,Y_{1+r,j'}>x)}{P(Y_{1,j}>x)^{1/\eta_{jj'}^{(r)}(\mb{Y})}}
\sim  \dst \frac{x^{-\alpha_j-\alpha_{j'}}+x^{-1-\alpha_j}+x^{-1-\alpha_{j'}}+x^{-1/\eta_{jj'}^{(r)}(\mb{X})}L(x)}
{(x^{-1}+x^{-\alpha_j})^{1/\eta_{jj'}^{(r)}(\mb{Y})}}
\end{array}
\end{eqnarray}
which implies the assertion. \fdem\\\\


In the EM4 case, the process $\{\mb{X}_n\}_{n\geq 1}$ is an M4 and thus $\eta_{jj'}^{(r)}(\mb{X})=1$ and $\alpha_j=\alpha$, for all $j=1,...,d$. Hence, we have
\begin{eqnarray}\nn
\eta_{jj'}^{(r)}(\mb{Y})=\left\{
\begin{array}{ll}
\dst \max\Big(\frac{1}{2},\alpha\Big)&,\,\alpha< 1\\
1&,\,\alpha\geq 1\,.
\end{array}\right.
\end{eqnarray}
already derived in Heffernan \emph{et al.} (\cite{zangh+} 2007; Result 10).\\

Observe that, in a pMAX process, $\alpha_j<1$ leads to tail independent random pairs $(Y_j,Y_{j'})$ for all $j'\in D$, whereas for $\alpha_j\geq 1$, random pairs $(Y_j,Y_{j'})$ will be tail dependent (independent) wether $(X_j,X_{j'})$ are tail dependent (independent). This is one important advantageous of the pMAX: it is  suitable for  multivariate time series data for which some pairs of components may present tail dependence while others may not.


\begin{ex}
Consider $\{\mb{X}_n\}_{n\geq 1}=\{(X_{n,1},X_{n,2},X_{n,3})\}_{n\geq 1}$ a $3$-dimensional sequence with $X_{n,1}=U_n$, $X_{n,2}=W_n$ and $X_{n,3}=W_{n+1}$, where sequences $\mb{U}=\{U_n\}_{n\geq 1}$ and  $\mb{W}=\{W_n\}_{n\geq 1}$ are such that:
\begin{itemize}
\item[(i)] $\mb{U}$ is i.i.d.\;with common d.f.\;$F_U$;
\item[(ii)] $F_U\in\mathcal{D}(G_1,\{a_n>0\},\{b_n\})$, i.e., $F_U$ belongs to domain of attraction of a non degenerate extreme values distribution $G_1$ with normalizing constants $\{a_n>0\}$ and $\{b_n\}$;
\item[(iii)]  $\mb{U}$ is independent of a chain $\mb{J}=\{J_n\}_{n\geq 1}$ of independent Bernoulli($1/2$) r.v.'s;
\item[(iii)]  $W_n=U_n\mathds{1}_{\{J_n=0\}}+U_{n+1}\mathds{1}_{\{J_n=1\}}$, $n\geq 1$.
\end{itemize}
Observe that $\mb{U}$ and $\mb{W}$ have the same common d.f.\;and thus normalized levels $u_n^{(\tau)}$ for $\mb{U}$ are also for $\mb{W}$. By hypothesis, $F_{X_{n,j}}^n(a_nx+b_n)\to G_1(x)$, $j=1,2,3$, and we have
$$
F_{\mb{X}_n}(x_1,x_2,x_3)=\frac{1}{4}\Big(2F_U(x_1\wedge x_2)F_U(x_3)+F_U(x_1)F_U(x_2\wedge x_3)+\prod_{j=1}^3F_U(x_j)\Big).
$$
Now observe that
\begin{eqnarray}\label{ex1}
\begin{array}{rl}
&n(1-F_{\mb{X}_n}(a_nx_1+b_n,a_nx_2+b_n,a_nx_3+b_n))
\vspace{0.35cm}\\
=&\frac{n}{4}\Big(4-2F_U(a_n(x_1\wedge x_2)+b_n)F_U(a_nx_3+b_n)-F_U(a_nx_1+b_n)F_U(a_n(x_2\wedge x_3)+b_n)\vspace{0.15cm}\\
&\hspace{0.5cm}-\prod_{j=1}^3F_U(a_nx_j+b_n)\Big)
\vspace{0.35cm}\\
=&\frac{n}{4}\Big(2\big(1-F_U(a_n(x_1\wedge x_2)+b_n)+F_U(a_n(x_1\wedge x_2)+b_n)(1-F_U(a_nx_3+b_n))\big)\vspace{0.15cm}\\
&\hspace{0.5cm}1-F_U(a_nx_1+b_n)+F_U(a_nx_1+b_n)(1-F_U(a_n(x_2\wedge x_3)+b_n))\vspace{0.15cm}\\
&\hspace{0.5cm}1-F_U(a_nx_1+b_n)+F_U(a_nx_1+b_n)(1-F_U(a_nx_2+b_n))\vspace{0.15cm}\\
&\hspace{0.5cm}+
F_U(a_nx_1+b_n)F_U(a_nx_2+b_n)(1-F_U(a_nx_3+b_n))\Big)\vspace{0.35cm}\\
\dst\mathop{\to}_{n\to\infty}&
\frac{1}{2}\big(-\log G_1(x_1\wedge x_2)-\log G_1(x_3)\big)\vspace{0.15cm}\\
&\hspace{0.5cm}+\frac{1}{4}\big(-\log G_1(x_1)-\log G_1( x_2\wedge x_3)\big)\vspace{0.15cm}\\
&\hspace{0.5cm}\frac{1}{4}\big(-\log G_1(x_1)-\log G_1(x_2)-\log G_1(x_3)\big).
\end{array}
\end{eqnarray}
Therefore,
\begin{eqnarray}\nn
\begin{array}{rl}
&F^n_{\mb{X}_n}(a_nx_1+b_n,a_nx_2+b_n,a_nx_3+b_n)
\vspace{0.35cm}\\
\dst\mathop{\to}_{n\to\infty}&
G(x_1,x_2,x_3)=G_1^{1/2}(x_1\wedge x_2)G_1^{1/4}( x_2\wedge x_3)G_1^{1/2}(x_1)G_1^{1/4}(x_2)G_1^{3/4}(x_3),
\end{array}
\end{eqnarray}
from which we derive the copula
\begin{eqnarray}\nn
\begin{array}{rl}
C_G(u_1,u_2,u_3)=(u_1\wedge u_2)^{1/2}( u_2\wedge u_3)^{1/4}u_1^{1/2}u_2^{1/4}u_3^{3/4}.
\end{array}
\end{eqnarray}
Next we compute the multivariate extremal index of $\{\mb{X}_n\}_{n\geq 1}$, which is a $2$-dependent sequence and thus satisfies a multivariate version of $D^{(3)}(\mb{u}^{(\boldsymbol{\tau})})$ condition of Chernick \emph{et al.} (\cite{chern}, 1991),
\begin{eqnarray}\nn
\begin{array}{rl}
\dst\lim_{n\to\infty}\sum_{j=4}^{[n/k_n]}P(\mb{X}_1\not\leq \mb{u}_n^{(\boldsymbol{\tau})},\mb{X}_2\leq \mb{u}_n^{(\boldsymbol{\tau})},\mb{X}_3\leq \mb{u}_n^{(\boldsymbol{\tau})},\mb{X}_j\not\leq \mb{u}_n^{(\boldsymbol{\tau})})=0,
\end{array}
\end{eqnarray}
for any sequence $\{k_n\}_{n\geq 1}$ such that, $k_n\to \infty$, $k_n\alpha_{n,l_n}\to 0$, $k_nl_n/n\to 0$, where $\{\alpha_{n,l}\}_{n\geq 1}$ is the sequence of strong-mixing coefficients and $\{l_n\}_{n\geq 1}$ is such that $\alpha_{n,l_n}\to 0$. Therefore, the multivariate extremal index is given by the following limit, whenever it exists:
\begin{eqnarray}\nn
\begin{array}{rl}
\dst\lim_{n\to\infty}\frac{P(\mb{X}_1\not\leq \mb{u}_n^{(\boldsymbol{\tau})},\mb{X}_2\leq \mb{u}_n^{(\boldsymbol{\tau})},\mb{X}_3\leq \mb{u}_n^{(\boldsymbol{\tau})})}{P(\mb{X}_1\not\leq \mb{u}_n^{(\boldsymbol{\tau})})}=\theta^{\mb{X}}(\tau_1,\tau_2,\tau_3).
\end{array}
\end{eqnarray}
Observe that the denominator can be derived from (\ref{ex1}), leading to
\begin{eqnarray}\nn
\begin{array}{rl}
\dst\lim_{n\to\infty}n P(\mb{X}_1\not\leq \mb{u}_n^{(\boldsymbol{\tau})})=\frac{1}{4}\big(2(\tau_1\vee\tau_2)
+(\tau_2\vee\tau_3)+2\tau_1+\tau_2+3\tau_3\big).
\end{array}
\end{eqnarray}
For the numerator, after some lengthy but simple calculations we obtain
\begin{eqnarray}\nn
\begin{array}{rl}
\dst\lim_{n\to\infty}n P(\mb{X}_1\not\leq \mb{u}_n^{(\boldsymbol{\tau})},\mb{X}_2\leq \mb{u}_n^{(\boldsymbol{\tau})},\mb{X}_3\leq \mb{u}_n^{(\boldsymbol{\tau})})=
\left\{
\begin{array}{ll}
4\tau_1 & ,\, \tau_1\geq \tau_2\geq\tau_3, \tau_1\geq \tau_3\geq\tau_2\\
\tau_1+3\tau_3 & ,\, \tau_3\geq \tau_1\geq\tau_2,\tau_3\geq \tau_2\geq\tau_1 \\
\tau_1+3\tau_2 & ,\, \tau_2\geq \tau_1\geq\tau_3, \tau_2\geq \tau_3\geq\tau_1,
\end{array}
\right.
\end{array}
\end{eqnarray}
and thus we have
\begin{eqnarray}\nn
\theta^{\mb{X}}(\tau_1,\tau_2,\tau_3)=
\left\{
\begin{array}{ll}
\frac{4\tau_1}{4\tau_1+2\tau_2+3\tau_3} & ,\, \tau_1\geq \tau_2\geq\tau_3\vspace{0.15cm}\\
\frac{\tau_1+3\tau_3}{4\tau_1+\tau_2+4\tau_3} & ,\, \tau_3\geq \tau_1\geq\tau_2\vspace{0.15cm}\\
\frac{4\tau_1}{4\tau_1+\tau_2+4\tau_3} & ,\, \tau_1\geq \tau_3\geq\tau_2\vspace{0.15cm}\\
\frac{\tau_1+3\tau_2}{2\tau_1+4\tau_2+3\tau_3} & ,\, \tau_2\geq \tau_1\geq\tau_3, \tau_2\geq \tau_3\geq\tau_1\vspace{0.15cm}\\
\frac{\tau_1+3\tau_3}{2\tau_1+3\tau_2+4\tau_3} & ,\, \tau_3\geq \tau_2\geq\tau_1.
\end{array}
\right.
\end{eqnarray}
as well as,
\begin{eqnarray}\nn
\begin{array}{cc}
\theta_{1,2}^{\mb{X}}(\tau_1,\tau_2)=
\left\{
\begin{array}{ll}
\frac{4\tau_1}{4\tau_1+2\tau_2} & ,\, \tau_1\geq \tau_2\vspace{0.15cm}\\
\frac{\tau_1+3\tau_2}{2\tau_1+4\tau_2} & ,\, \tau_2\geq \tau_1
\end{array}
\right.,
&
\theta_{1,3}^{\mb{X}}(\tau_1,\tau_3)=
\left\{
\begin{array}{ll}
\frac{\tau_1+3\tau_3}{4\tau_1+4\tau_3} & ,\, \tau_3\geq \tau_1\vspace{0.15cm}\\
\frac{4\tau_1}{4\tau_1+4\tau_3} & ,\, \tau_1\geq \tau_3
\end{array}
\right.,
\end{array}
\end{eqnarray}
\begin{eqnarray}\nn
\begin{array}{c}
\theta_{2,3}^{\mb{X}}(\tau_2,\tau_3)=
\left\{
\begin{array}{ll}
\frac{3\tau_3}{3\tau_2+4\tau_3} & ,\, \tau_3\geq \tau_2\vspace{0.15cm}\\
\frac{3\tau_2}{4\tau_2+3\tau_3} & ,\, \tau_2\geq \tau_3,
\end{array}
\right.
\end{array}
\end{eqnarray}
and univariate extremal indexes, $\theta_1^{\mb{X}}=1$ and $\theta_2^{\mb{X}}=\theta_3^{\mb{X}}=3/4$.

Applying the lag-$r$ TDC definition in (\ref{lambdar}), simple calculations lead us to the results summarized in Table \ref{tab1}.
\begin{table}
\caption{The lag-$r$ TDC coefficients of sequence $\{\mb{X}_n\}_{n\geq 1}$.}\label{tab1}
\begin{center}
\begin{tabular}{c|ccc}
$\lambda_{jj'}^{(r)}(\mb{X})$&j'=1 & j'=2 &j'=3\\
\hline
j=1 & $0$ if $r\geq 1$ & $1/2$ if $r=0$& $0$  if $r\geq 0$\\
 & & $0$  if $r\geq 1$& \\
 \hline
j=2 & $1/2$ if $r=1$ &$1/4$ if $r=1$ & $1/4$ if $r=0$  \\
& $0$ if $r\geq 2$ & $0$ if $r\geq 2$ &  $0$ if $r\geq 1$  \\
\hline
j=3 & $1/2$ if $r=1,2$ & $1/2$ if $r=1$ & $1/4$ if $r=1$ \\
 &  $0$ if $r\geq 3$ &  $1/4$ if $r=2$   &  $0$ if $r\geq 2$ \\
 &  &  $0$ if $r\geq 3$   & \\
 \hline
\end{tabular}
\end{center}
\end{table}
Obviously $\eta_{jj'}^{(r)}(\mb{X})=1$ whenever $\lambda_{jj'}^{(r)}(\mb{X})$ is positive and, since we obtain $P(X_{1+r,j'}>x,X_{1,j}>x)=P(X_{1,j}>x)^2$ in all cases of a null TDC, we have $\eta_{jj'}^{(r)}(\mb{X})=1/2$ in those cases.

Now if we consider a pMAX process
\begin{eqnarray}\nn
\begin{array}{lc}
\{\mb{Y}_n=(Y_{n,1},Y_{n,2},Y_{n,3})\}_{n\geq 1}=\{(X_{n,1}\vee Z_{n,1}^{1/\alpha_1},(X_{n,2}\vee Z_{n,2}^{1/\alpha_2},X_{n,3}\vee Z_{n,3}^{1/\alpha_3})\}_{n\geq 1}
\end{array}
\end{eqnarray}
with $\alpha_1=3/2$, $\alpha_2=1$ and $\alpha_3=2/3$, then by applying the Proposition \ref{ppmaxei}, we have $\theta_2^{\mb{Y}}=\theta_2^{\mb{X}}=3/4$, $\theta_1^{\mb{Y}}=\theta_3^{\mb{Y}}=\theta_{1,3}^{\mb{Y}}(\tau_1,\tau_3)=1$, $\theta_{1,2}^{\mb{Y}}(\tau_1,\tau_2)=\theta_{1,2}^{\mb{X}}(\tau_1,\tau_2)$, \begin{eqnarray}\nn
\begin{array}{ccc}
\theta_{2,3}^{\mb{Y}}(\tau_2,\tau_3)=\frac{\frac{3}{4}\tau_2+\tau_3}{\tau_2+\tau_3}
& \textrm{ and }
& \theta^{\mb{Y}}(\tau_1,\tau_2,\tau_3)=\frac{\theta_{1,2}^{\mb{X}}(\tau_1,\tau_2)
(\tau_1\vee\tau_2)+\tau_1+\tau_2+2\tau_3}
{(\tau_1\vee\tau_2)+\tau_1+\tau_2+2\tau_3}
=
\left\{
\begin{array}{ll}
\frac{4\tau_1^2+\tau_2^2+3\tau_1\tau_2+4\tau_1\tau_3+2\tau_2\tau_3}
{4\tau_1^2+\tau_2^2+4\tau_1\tau_2+4\tau_1\tau_3+2\tau_2\tau_3}
& ,\, \tau_1\geq \tau_2\vspace{0.15cm}\\ \frac{2\tau_1^2+7\tau_2^2+7\tau_1\tau_2+4\tau_1\tau_3+8\tau_2\tau_3}
{2\tau_1^2+8\tau_2^2+8\tau_1\tau_2+4\tau_1\tau_3+8\tau_2\tau_3} & ,\, \tau_2\geq \tau_1.
\end{array}
\right.
\end{array}
\end{eqnarray}
The pMAX lag-$r$ TDC's are $\lambda_{1j'}^{(r)}(\mb{Y})=\lambda_{1j'}^{(r)}(\mb{X})$, $\lambda_{2j'}^{(r)}(\mb{Y})=\frac{1}{2}\lambda_{2j'}^{(r)}(\mb{X})$ and $\lambda_{3j'}^{(r)}(\mb{Y})=0$, for all $r\in \mathbb{N}_0$ and $j'=1,2,3$. In what concerns the lag-$r$ Ledford and Tawn coefficients, we have $\eta_{j1}^{(r)}(\mb{Y})=\eta_{j2}^{(r)}(\mb{Y})=1/2$ if $j=1,2$ and $r\geq j$, $\eta_{j3}^{(r)}(\mb{Y})=3/5$ if $j=1,2$ and $r\geq j-1$, $\eta_{3j'}^{(r)}(\mb{Y})=2/5$ if $j=1,2$ and $r\geq 3$ and $\eta_{33}^{(r)}(\mb{Y})=1/2$ if $r\geq 2$.
\end{ex}

\section{Estimation: some notes}\label{sestim}


In this section we shall give some guide marks on the estimation of the marginal parameters $\alpha_j$, $j\in D$, and the tail dependence coefficients within a pMAX process. For the estimation of the multivariate extremal index, we refer the manuscript of Smith and Weissman (\cite{smith+weiss}, 2006) and the paper of Robert (\cite{rob}, 2008).\\

We start by remarking that $\alpha_j$, $j\in D$, are the key parameters of the pMAX process dictating the type of dependence and thus allowing to sometimes deduce the value of some coefficients. For instance, an estimate of $\alpha_j$ less than one means that the random pairs $(Y_{n,j},Y_{n+r,j'})$, for all $n,r\in \en$ and $j'\in D$, are tail independent leading to a respective null lag-$r$ TDC, $\lambda_{jj'}^{(r)}(\mb{Y})$. In addition,  estimates of $\alpha_j$ less than one for all $j\in D$ will imply an unit multivariate extremal index. So, as a first step, we suggest to marginally apply a test of tail independence (see e.g, Zhang \cite{zhang} 2008 and references therein).

In order to estimate the value of $\alpha_j$, observe that it corresponds to the tail index of the marginal process $\{Y_{n,j}\}_{n\geq 1}$ whenever $\alpha_j\leq 1$. There are several known  tail index estimators in literature
such as, Hill, Pickands, maximum likelihood estimator,
moments, generalized
weighted moments, among
others. Their properties have been derived in an i.i.d. framework,
but there are some studies considering a stationary context (see,
for instance, Rootz\'{e}n {et al.}  \cite{root+} 1990, Resnick and
St\v{a}ric\v{a} \cite{res+sta1,res+sta2} 1995/1998, Drees
\cite{drees1} 2003). The Hill estimator is the most used in the context of a Fréchet domain of attraction. In particular, under a strong-mixing dependence structure, it is consistent and asymptotically
normal (Rootz\'{e}n {et al.} \cite{root+}).

For a given sequence $\{\xi_n\}_{n\geq 1}$, a strong-mixing condition means that $\xi_n$ and $\xi_{n+s}$ are approximately independent for large $s$. More formally, it will hold whenever $$|P(A \cap B) -P(A) P(B)|\leq a(s)\dst\mathop{\to}_{s\to\infty} 0$$
for $A\in \mathcal{F}( \xi_1,...,\xi_k)$, $B\in \mathcal{F}( \xi_{k+s},\xi_{k+s+1},...)$ and $k,s\geq 1$, where $ \mathcal{F}( \xi_{1},...,\xi_{k},...)$ denotes the $\sigma$-algebra generated by $\{\xi_{1},...,\xi_{k},...\}$.
This is a mild condition to assume as many stationary processes, including M4 and EM4, do satisfy strong-mixing. In our case it will be easily stated for the marginal process $\{Y_{n,j}\}_{n\geq 1}$ by just assuming that it holds for the underlying process $\{X_{n,j}\}_{n\geq 1}$ (Bradley (\cite{brad}, 2005; Theorem 5.2).

We can also find in literature estimation procedures of the TDC and of the Ledford and Tawn coefficient. However, the properties of consistency and asymptotic normality are  mainly derived under an i.i.d.\, assumption, i.e., considering i.i.d.\;random pairs, which is not an interesting case when we intend to estimate the respective lag-$r$ versions, $\lambda_{jj'}^{(r)}(\mb{Y})$ and $\eta_{jj'}^{(r)}(\mb{Y})$.
A brief note on the estimation of the TDC, when we drop the independence assumption, can be seen in Ferreira and Ferreira (\cite{hf+mf}, 2012). In what concerns coefficient $\eta_{jj'}^{(r)}(\mb{Y})$ defined in (\ref{etar}), observe that it can be estimated as the tail index of $\min(Y_{1,j},Y_{1+r,j'})$. By assuming that the underlying sequence $\{\mb{X}_n\}_{n\geq 1}$ satisfies the strong-mixing property we also derive a strong-mixing structure for sequence $\{\min(Y_{n,j},Y_{n+r,j'})\}_{n\geq 1}$ (see Lemmas \ref{lem1} and \ref{lem2} in Appendix). Now the considerations above concerning the Hill estimator will also hold in this case.

\appendix
\section{Appendix}

\begin{lem}\label{lem1}
Let $\{\mb{Y}_n\}_{n\geq 1}$ be a $d$-dimensional sequence satisfying the strong-mixing condition with a sequence of dependence coefficients $\{a(s)\}_{s\geq 1}$. Then,
\begin{enumerate}
\item[(i)]$\{(Y_{n,j},Y_{n+r,j'})\}_{n\geq 1}$ satisfies the strong-mixing condition with any sequence of dependence coefficients $\{b(s)\}_{s\geq 1}$ such that $b(s)=a(s-r)$, $s>r$.
\item[(ii)] $\{(Y_{n,j}\wedge Y_{n+r,j'})\}_{n\geq 1}$ satisfies the strong-mixing condition with $\{b(s)\}_{s\geq 1}$ as in (i).
\end{enumerate}
\end{lem}
\dem
\begin{enumerate}
\item[(i)] By hypothesis, $\forall A\in \mathcal{F}(\mb{Y}_1,...,\mb{Y}_p),\, B\in \mathcal{F}(\mb{Y}_{p+s},\mb{Y}_{p+s+1},...)$, $s\geq 1$, we have
$$
|P(A\cap B)-P(A)P(B)|\leq a(s)\dst\mathop{\to}_{s\to\infty} 0.
$$
Consider $r\geq 1$ fixed, $C\in \mathcal{F}((Y_{1,j},Y_{1+r,j'}),...,(Y_{p,j},Y_{p+r,j'})),\, D\in \mathcal{F}((Y_{p+s,j},Y_{p+r+s,j'}),...)$ and $s>r$. We have
$
C\in \mathcal{F}(\mb{Y}_1,...,\mb{Y}_{p+r}),\, D\in \mathcal{F}(\mb{Y}_{p+s},\mb{Y}_{p+s+1},...)
$
and hence
$$
|P(C\cap D)-P(C)P(D)|\leq a(s-r)\dst\mathop{\to}_{s\to\infty} 0.\,\square
$$
\item[(ii)] Just observe that
$$
\mathcal{F}(Y_{i,j}\wedge Y_{i+r,j'},\,i\in I)\subset \mathcal{F}((Y_{i,j}, Y_{i+r,j'}),\,i\in I)
$$
for any family of indexes $I$. \fdem
\end{enumerate}

\begin{lem}\label{lem2}
Let $\{\mb{X}_n\}_{n\geq 1}$ be a $d$-dimensional sequence satisfying the strong-mixing condition and $\{\mb{Z}_n\}_{n\geq 1}$ a $d$-dimensional i.i.d.\,sequence independent of the previous. Then $\{\mb{Y}_n=\mb{X}_n\vee \mb{Z}_n\}_{n\geq 1}$ satisfies the strong-mixing condition.
\end{lem}
\dem
Consider sequence $\{\mb{U}_n=(X_{n,1},...,X_{n,d}, Z_{n,1},...,Z_{n,d})\}_{n\geq 1}$ and let $A\in \mathcal{F}(\mb{U}_1,...,\mb{U}_p),\, B\in \mathcal{F}(\mb{U}_{p+s},\mb{U}_{p+s+1},...)$. Since we can write
$$
|P(A\cap B)-P(A)P(B)|=|P(A'\cap B')P(A''\cap B'')-P(A')P(A'')P(B')P(B'')|,
$$
with $A'\in \mathcal{F}(\mb{X}_1,...,\mb{X}_p),\, B'\in \mathcal{F}(\mb{X}_{p+s},\mb{X}_{p+s+1},...),\,A''\in \mathcal{F}(\mb{Z}_1,...,\mb{Z}_p),\, B''\in \mathcal{F}(\mb{Z}_{p+s},\mb{Z}_{p+s+1},...)$, we have that
$$
\begin{array}{rl}
|P(A\cap B)-P(A)P(B)|\leq& |P(A'\cap B')-P(A')P(B')|P(A''\cap B'')\vspace{0.15cm}
\\
&+|P(A''\cap B'')-P(A'')P(B'')|P(A')P(B')\vspace{0.35cm}
\\
\leq &a_{\mb{X}_n}(s)\dst\mathop{\to}_{s\to\infty} 0.
\end{array}$$
Therefore sequence $\{\mb{U}_n\}_{n\geq 1}$ satisfies the strong-mixing property. Now the result follows from the same argument used in Lemma \ref{lem1}.(ii) for the maximum operator ($\vee$).

\end{document}